\documentclass{amsart}
\usepackage{amsmath,amssymb,amscd,amsthm}
\usepackage{enumerate}
\usepackage[dvips]{graphicx,color}

\newtheorem{theorem}{Theorem}[section]
\newtheorem{lemma}[theorem]{Lemma}
\newtheorem{proposition}[theorem]{Proposition}
\newtheorem{assertion}[theorem]{Assertion}

\newtheorem{mtheorem}{Theorem}

\newtheorem{mcorollary}[mtheorem]{Corollary}

\numberwithin{figure}{section}
\numberwithin{equation}{section}
\theoremstyle{definition}
\newtheorem{definition}[theorem]{Definition}
\theoremstyle{remark}
\newtheorem{remark}[theorem]{Remark}

\begin{document}

\title[Heterodimensional Tangencies and Strange Attractors]{Heterodimensional Tangencies on cycles Leading To Strange Attractors}

\author{Shin Kiriki}
\address{Department of Mathematics, Kyoto University of Education, 
1 Fukakusa-Fujinomori, Fushimi-ku, Kyoto, 612-8522, Japan}
\email{skiriki@kyokyo-u.ac.jp}

\author{Yusuke Nishizawa}
\address{Department of Mathematics and Information Sciences, Tokyo Metropolitan
University, Minami-Ohsawa 1-1, Hachioji, Tokyo 192-0397, Japan}
\email{nishizawa-yusuke@ed.tmu.ac.jp}

\author{Teruhiko Soma}
\address{Department of Mathematics and Information Sciences, Tokyo Metropolitan
University, Minami-Ohsawa 1-1, Hachioji, Tokyo 192-0397, Japan}
\email{tsoma@tmu.ac.jp}

\date{\today,\quad to appear in \emph{Discrete Conti. Dynam. Sys.} } 


\maketitle

\begin{abstract}
In this paper, we study a two-parameter family $\{\varphi_{\mu,\nu}\}$ of 3-dimensional 
diffeomorphisms which have a  bifurcation 
induced by simultaneous generation of a heterodimensional cycle and a heterodimensional tangency 
associated to two saddle points. 
We show that such a codimension-$2$ bifurcation generates a quadratic homoclinic tangency 
associated to one of the saddle continuations which unfolds generically with respect to some  
one-parameter subfamily of $\{\varphi_{\mu,\nu}\}$.
Moreover, from this result together with some well-known  facts, 
we detect some nonhyperbolic phenomena (i.e.,\ the existence of 
nonhyperbolic strange attractors 
and the $C^{2}$ robust tangencies) arbitrarily close to the codimension-$2$ bifurcation.
\end{abstract}

\section{Introduction}\label{introduction}
When diffeomorphisms act on manifolds of dimension greater than or equal to three,
it is well known that nonhyperbolic phenomena 
are caused by the existence of heteroclinic cycles containing two saddle points with 
different indexes, called \textit{heterodimensional cycles}, as well as that of
homoclinic tangencies \cite{N79, R83, GTS93, PV94, R94, GST96}.
Heterodimensional cycles presented a new mechanism in dynamics, which has been studied D\'iaz et al.\ \cite{DR92, 
D95_1,DR97,BD99,DR01,BDP03,BD08}.
As is suggested in \cite{P05}, these two classes of nonhyperbolic diffeomorphisms seem to occupy 
a large part in the complement of the hyperbolic diffeomorphisms. 
In this paper, we study 
$3$-dimensional diffeomorphisms which have heterodimensional cycles and tangencies of certain type 
simultaneously. 

Let $\varphi$ be a diffeomorphism on a $3$-dimensional smooth manifold $M$ which has 
two saddle fixed points $p$ and $q$ satisfying $\mathrm{index}(q)= \mathrm{index}(p)+1$, 
where $\mathrm{index}(\cdot)$ denotes the dimension of the unstable manifold of the saddle point.
A heteroclinic point $r$ of the stable manifold $W^s(p)$ and the unstable 
manifold $W^u(q)$ is called a \emph{heterodimensional tangency} of $W^s(p)$ and $W^u(q)$ 
if $r$ satisfies   
\begin{itemize}
\item$T_r W^s(p)+T_r W^u(q)\not=T_rM$;
\item $\dim(T_r W^s(p)) + \dim(T_r W^u(q))>\dim(M)$,
\end{itemize}
where $\dim(\cdot)$ denotes the dimension of the space.

Our main theorem in this paper is as follows.
We will present some definitions  and  generic conditions used in the statement of Theorem \ref{thm_A} in Section 
\ref{Generic_conditions}.

\begin{mtheorem}\label{thm_A}
Let $M$ be a $3$-dimensional $C^2$ manifold, and let
$\{\varphi_{\mu, \nu}\}$ be a two-parameter family of  
$C^2$ diffeomorphisms $\varphi_{\mu,\nu}:M\rightarrow M$ which $C^2$ depends on $(\mu,\nu)$ and has continuations 
of saddle fixed points 
$p_{\mu,\nu}$ and $q_{\mu,\nu}$ with $\mathrm{index}(p_{\mu,\nu})=1$ and $\mathrm{index}(q_{\mu,\nu})=2$.
Suppose that the following conditions hold.
\begin{itemize}
\item
Each $\varphi_{\mu,\nu}$ is locally $C^2$ linearizable in a small neighborhood $N(q_{\mu,\nu})$ 
of $q_{\mu,\nu}$.
\item
$\varphi=\varphi_{0,0}$ has a heterodimensional cycle containing the fixed points $p=p_{0,0}$, $q=q_{0,0}$, 
a nondegenerate heterodimensional tangency $r$, a quasi-transverse intersection $s\in W^s(q)\cap W^u(p)$.
\item
$\{\varphi_{\mu,\nu}\}$ satisfies the generic conditions (C1)-(C4) given in Section 
\ref{Generic_conditions}.
\end{itemize}
Then, for a sufficiently small $\varepsilon>0$ and any $\mu$ in either $(0,\varepsilon)$ or $(-\varepsilon,0)$, 
there exist infinitely many $\nu$ such that 
$\varphi_{\mu, \nu}$ has a quadratic homoclinic tangency associated to $p_{\mu, \nu}$ which unfolds generically 
with respect to the $\nu$-parameter family $\{ \varphi_{\mu(\mathrm{fixed}),\nu} \}$.
\end{mtheorem}

\begin{figure}[hbtp]
\centering
\scalebox{0.55}{\includegraphics[clip]{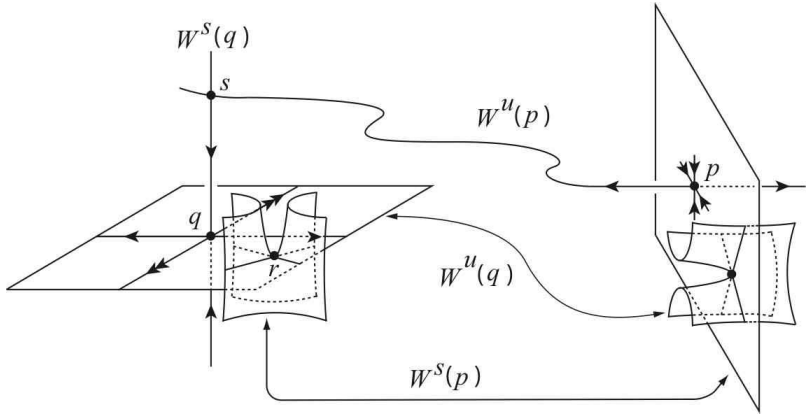}}
\caption{A heterodimensional cycle associated with the pair of saddles $p$ and $q$, which contains 
 the heterodimensional tangency $r$ 
and the quasi-transverse intersection $s$.}
\label{fg_1}
\end{figure}

\begin{remark}
\begin{enumerate}
\item The generic conditions (C2) and (C3) imply the setting such that 
a heterodimensional tangency and a quasi-transverse intersection
unfold generically with respect to the parameters $\mu$ and $\nu$, respectively.
The point to notice is that the newly detected  homoclinic tangency in Theorem \ref{thm_A} can be 
also controlled by these parameters.
It is not hard to get a one parameter family in the infinite dimensional space $\mathrm{Diff}^2(M)$ with respect to which the homoclinic tangency unfolds generically.
However, in our proof, we need to show that the tangency given in Subsection \ref{Existence} unfolds 
generically with respect to the $\nu$-parameter family $\{\varphi_{\mu_0,\nu}\}$ in the two-dimensional subspace $\{\varphi_{\mu,\nu}\}$ of $\mathrm{Diff}^2(M)$.

\item Theorem \ref{thm_A} holds for a homoclinic tangency 
associated to $q_{\mu,\nu}$ instead of $p_{\mu,\nu}$ if we replace 
the conditions in (C1) and (C4) by appropriate ones.
\item  We need the $C^2$ smoothness in the local linearization around $q_{\mu,\nu}$ to estimate the curvatures 
of $W^s(p_{\mu,\nu})$ and $W^u(p_{\mu,\nu})$ at the tangency in Section \ref{some_lemma}.
\end{enumerate}
\end{remark}
Figure \ref{fg_1} illustrates an example of 
heterodimensional cycles containing a nondegenerate heterodimensional tangency of hyperbolic type, see Definition \ref{d_generic_unfolding}-(3).

The conclusion obtained from Theorem \ref{thm_A} reminds us of prior works associated with homoclinic tangencies.
The one is related to strange attractors and the other $C^{2}$ robust tangencies. 

First, let us discuss the former one. Viana showed the following theorem.
\begin{theorem}[Viana \cite{V93}]\label{Viana}
For a generic subset of one-parameter families $\{ {\varphi}_{\mu}\}$ of $C^\infty$ diffeomorphisms 
on any manifolds of the dimension greater than or equal to two that unfolds a homoclinic tangency 
at parameter value $\mu=0$ associated to a sectionally dissipative saddle periodic point, there 
is a subset $S$ of $\mathbb{R}$ such that 
\begin{itemize}
\item $S \cap (-\epsilon, \epsilon)$ has a positive Lebesgue measure for every $\epsilon >0$, 
\item for all $\mu \in S$, ${\varphi}_{\mu}$ exhibits nonhyperbolic strange attractors in a 
$\mu$-dependent neighborhood of the orbit of tangency.  
\end{itemize}
\end{theorem}

Leal \cite{L08} extended this result and showed the existence of infinitely many strange attractors.
A saddle periodic point is said to be \emph{sectionally dissipative} if  
the product of any two eigenvalues of the derivative at the  point has norm less than one. 
Also, $\Lambda$ is a \emph{strange attractor} of $\varphi$ if $\Lambda$ is a compact, $\varphi$-invariant, transitive set and  
the basin $W^{s}(\Lambda)$ has a nonempty interior and 
there exists $z_{1}\in \Lambda$ such that $\{ {\varphi^{n}(z_{1})}:n \geq 0\}$ is dense in $\Lambda$ and 
$|| d{\varphi}^{n}(z_{1}) || \geq e^{cn}$ for all $n\geq 0$ and some $c>0$.
Note, Viana \cite{V93} assumed that for simplicity ${\varphi}_{\mu}$ is $C^{3}$ linearizable in a neighborhood of the saddle point for $\mu$ is sufficiently close to $0$.
Combining these extra conditions and our result for the cycle containing the heterodimensional tangency
imply the following corollary.  

\begin{mcorollary}\label{cor_B}
Let $\{\varphi_{\mu, \nu}\}$ be the two-parameter family of $3$-dimensional $C^{\infty}$ diffeomorphisms 
having the heterodimensional cycle with the nondegenerate heterodimensional tangency for $(\mu,\nu)=(0,0)$ 
given in Theorem \ref{thm_A}. 
Suppose that each $\varphi_{\mu,\nu}$ is locally $C^{3}$ linearizable in a small neighborhood of $p_{\mu,\nu}$ 
and $\varphi_{0,0}$ is sectionally dissipative at $p_{0,0}$.
Then there exists a positive Lebesgue measure subset 
$\mathcal{A}$ of the  $\mu\nu$-plane arbitrarily near $(0,0)$ such that, for any $(\mu, \nu)\in \mathcal{A}$, 
$\varphi_{\mu, \nu}$ exhibits nonhyperbolic strange attractors.
\end{mcorollary}

Next, we discuss $C^{2}$ robust homoclinic tangencies derived from the heterodimensional cycle having a 
nondegenerate  tangency. 
A homoclinic tangency of a diffeomorphism $\varphi$ associated with a hyperbolic set $\Gamma$ 
is $C^{r}$ \emph{robust} if there is a $C^{r}$ neighborhood $\mathcal{U}$ of  $\varphi$ 
such that every diffeomorphism $\psi\in \mathcal{U}$ has a homoclinic tangency associated with the continuations of $\Gamma$ for $\psi$. 
Newhouse \cite{N79} showed, in the $C^{2}$ topology, 
a homoclinic tangency of surface diffeomorphisms generates $C^{2}$ robust homoclinic tangencies.  
This property yields the so-called $C^{2}$ Newhouse phenomenon:  
there is a non-empty open set of $C^{2}$ diffeomorphisms
and its residual subset such that every diffeomorphism in the subset has infinitely many sinks.
The result is extended by Palis-Viana \cite{PV94} to  the higher dimensional case.
Palis and Viana proved the result under the sectional dissipativeness and linearizing conditions as in \cite{V93}. 
Moreover, Romero \cite{R94} proved the following theorem without these conditions.

\begin{theorem}[Romero \cite{R94}]\label{Romero}
Let $\varphi$ be a $C^2$ diffeomorphism on a manifold $M$, $\dim(M)\geq 3$, 
having a homoclinic tangency associated with 
a saddle periodic point  of $\varphi$ whose index is greater or equal to $2$.
Then there are diffeomorphisms arbitrarily $C^2$ close to $\varphi$ having robust
homoclinic tangencies. 
\end{theorem}

Therefore, we have the following corollary.

\begin{mcorollary}\label{cor_C}
Let $\{\varphi_{\mu, \nu}\}$ be the two-parameter family of $3$-dimensional $C^2$ diffeomorphisms  
having the heterodimensional cycle with the nondegenerate heterodimensional tangency
for $(\mu,\nu)=(0,0)$ given in Theorem \ref{thm_A}.   
Then, for a sufficiently small $\varepsilon>0$ and any $\mu$ in either $(0,\varepsilon)$ or $(-\varepsilon,0)$, 
there are infinitely many $\nu$ such that 
every $C^2$  neighborhood of $\varphi_{\mu,\nu}$ contains 
a diffeomorphism having a $C^2$ robust homoclinic tangency.
\end{mcorollary}

\begin{remark}\label{rem_0}
If we add a \emph{weak} dissipative condition (see in \cite[Theorem A (1.1)]{R94})
to the assumptions of Corollary \ref{cor_C} then, the $C^{2}$ Newhouse phenomenon is obtained from our settings.
\end{remark}

\begin{remark}\label{rem_1}
In the $C^1$ case,    
D\'iaz et al.\ \cite{DNP06} have already proved that  the unfolding of heterodimensional tangencies leads to 
non-dominated dynamics and therefore (by results of \cite{BD99}  and \cite{BDP03}) to the $C^{1}$ Newhouse  
phenomenon (see also \cite{As} for a different approach to the phenomenon). 
On the other hand, the theory of strange attractors has not been so far developed in the $C^l$ category $(l=1,2)$.
\end{remark}

\noindent
{\bf Outline of the proof of Theorem \ref{thm_A}}:
We will finish Introduction by presenting a sketch of the proof of Theorem \ref{thm_A}.

In Section \ref{Generic_conditions}, we give definitions and the generic conditions (C1)-(C4) used in Theorem \ref{thm_A}.
Especially, a nondegenerate heterodimensional tangency, which is one of main ingredients of this paper, is introduced explicitly there.
Such tangencies are classified into the elliptic and hyperbolic types, see  Definition \ref{d_generic_unfolding}-(3).

\begin{figure}[hbtp]
\centering
\scalebox{0.4}{\includegraphics[clip]{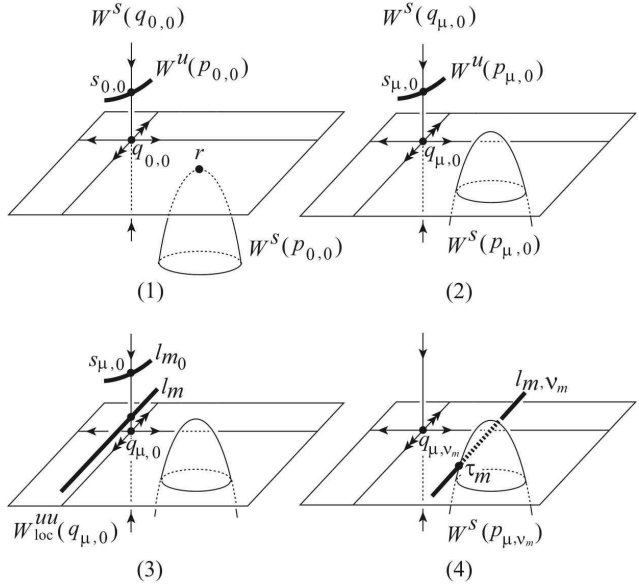}}
\caption{(1) The situation when $(\mu,\nu)=(0,0)$.
(2) $s_{\mu,0}$ is a quasi-transverse intersection  with respect to the new parameter.
(3) $l_{m}$ converges $W^{uu}_{\mathrm{loc}}(q_{{\mu}, 0})$ in $C^{2}$ topology.
(4) $l_{m,\nu_m}$ is a continuation of $l_m$ for $\nu=\nu_m$.}
\label{fg_9}
\end{figure}

In Section \ref{some_lemma}, we prove three lemmas.
Lemma \ref{new_parameters} shows the existence of the re-parametrization of $\{\varphi_{\mu, \nu}\}$ such that, 
for any $(\mu,0)$ in the new parameter, $W^s(q_{\mu,0})$ and $W^u(p_{u,0})$ still have a quasi-transverse intersection $s_{\mu,0}$ 
which unfolds generically with respect to the $\nu$-parameter but the tangency $r$ is annihilated.
See Fig.\ref{fg_9} (1) and (2).
Lemma \ref{inclination} applies the $C^{2}$ inclination lemma to a shorter curve $l_{m_0}$ in $W^u(p_{\mu,0})$ passing 
through $s_{\mu,0}$ so that some curve $l_m$ in $\varphi_{\mu,0}^m(l_{m_0})$ containing $\varphi_{\mu,0}^m(s_{\mu,0})$ $C^2$ converges to $W^{uu}_{\mathrm{loc}}(q_{\mu,0})$ as $m\rightarrow \infty$. 
See Fig.\ref{fg_9} (3).
Lemma \ref{quadratic proof} explains a connection between quadratic tangencies and curvatures, which is 
used to show that the homoclinic tangencies obtained in Section \ref{Proof} are quadratic.

Assertion \ref{elliptic} and Assertion \ref{hyperbolic} in Section \ref{Proof} show that the generic unfolding of heterodimensional tangencies introduces 
the existence of a quadratic homoclinic tangency $\tau_m$ as illustrated in Fig.\ \ref{fg_9}-(4).
Finally, we prove in Proposition \ref{prop2} that 
the tangency $\tau_m$ unfolds generically with respect to the ${\nu}$-parameter.
These results assure the proof of Theorem \ref{thm_A}.

\section{Definitions and generic conditions}\label{Generic_conditions}

In this section, we present some definitions needed in later sections and generic conditions adopted 
as hypotheses in Theorem \ref{thm_A}.

\subsection{Definitions}

\begin{definition}\label{d_generic_unfolding}
Suppose that $M$ is a 3-dimensional $C^2$ manifold.
Let $\{l_\nu\}_{\nu\in J}$,  $\{m_\nu\}_{\nu\in J}$ be  $C^2$ families of regular curves in $M$, and let 
$\{S_\nu\}_{\nu\in J}$, $\{Y_\nu\}_{\nu\in J}$ be $C^2$ families of regular surfaces in $M$, where 
$J$ is an open interval.

(1)
Suppose that $l_{\nu_0}$ and $m_{\nu_0}$ intersect at a point $s$ for some $\nu_0\in J$ and some open 
neighborhood $U$ of $s$ in $M$ has a $C^2$ change of  coordinates with respect to which $m_\nu=\{(0,0,z)\in U\}$ 
for any $\nu\in J$ near $\nu_0$.
We say that $s$ is a \emph{quasi-transverse intersection} of $l_{\nu_0}$ and $m_{\nu_0}$ if
$$\dim(T_{s}(l_{\nu_0}) + T_{s}(m_{\nu_0}))=2.$$
Moreover, $s$ \emph{unfolds generically} at $\nu = \nu_0$ with respect to $\{l_\nu\}_{\nu\in J}$, $\{m_\nu\}_{\nu\in J}$ if there exists 
a $C^2$ map $s:J\rightarrow M$ with $s(\nu) = s_{\nu} \in l_{\nu}$ for any $\nu \in J$ and $s({\nu}_{0})=s$ 
and a $C^2$ function $d:J\rightarrow \mathbb{R}^+$ with 
$d(\nu_0)\neq 0$ such that 
\begin{equation}\label{quasi_transverse}
T_{s} M=T_{s}( l_{\nu_{0}})\oplus N \oplus T_{s} (m_{\nu_{0}})\quad\mbox{and}\quad\mathrm{dist}(s_{\nu}, m_{\nu}) = 
\vert \nu-\nu_0\vert d(\nu)
\end{equation}
for any $\nu$ near $\nu_0$, 
where $N$ is the one-dimensional space spanned by the non-zero tangent vector 
$(d s_{\nu}/d \nu) |_{\nu=\nu_0}$.
This property corresponds to the conditions (GU1)--(GU3) in \cite[\S2.2.1]{DR01}.

(2)
Suppose that $l_{\nu_0}$ and $S_{\nu_0}$ intersect at a point $\tau$ for some $\nu_0\in J$.
We say that $\tau$ is a \emph{quadratic tangency} (or a \emph{contact of order} 
$1$) of $l_{\nu_0}$ and $S_{\nu_0}$ if there exists some $C^2$ change of coordinates on an open neighborhood $U(\tau)$ of $\tau$ 
such that $\tau= (0,0,0)$, $S_\nu=\{(x,y,z)\in U(\tau)\,;\, z = 0\}$ 
and $l_\nu$ has a regular parametrization 
$l(\nu,t)=(x(\nu,t),y(\nu,t),z(\nu,t))$ with $l(\nu_0,0)=(0,0,0)$ and 
$$\frac{\partial z}{\partial t}(\nu_0,0)= 0\quad\mbox{and}\quad \displaystyle \frac{\partial^2 z}{\partial t^2}
(\nu_0,0)\neq 0,$$

The tangency $\tau$ is said to \emph{unfold generically} at $\nu = \nu_0$ with respect to 
$\{l_\nu\}_{\nu\in J}$ and $\{S_\nu\}_{\nu\in J}$ if
$$\frac{\partial z}{\partial \nu}(\nu_0,0)\neq 0.$$

(3)
Suppose that $S_{\nu_0}$ and $Y_{\nu_0}$ intersect at a point $r$ for some $\nu_0\in J$.
We say that $r$ is a \emph{nondegenerate heterodimensional tangency} of $S_{\nu_0}$ and $Y_{\nu_0}$ 
if  there exists a $C^2$ coordinate on an open set $U$ in $M$ containing $r$  
with $r= (u_0,v_0,0)$ 
for some $u_0,v_0\in\mathbb{R}$, $S_\nu=\{(x,y,z)\in U\,;\, z = 0\}$ and such that $Y_\nu$ has a 
parametrization $(x,y,f_\nu(x,y))$ the third entry $f_\nu(x,y)=f(\nu,x,y)$ of which  is a $C^2$ function satisfying 
\begin{equation}\label{nondegenerate_condition}
f_{\nu_0}(u_0,v_0)=0,\quad
\frac{\partial f_{\nu_0}}{\partial x}(u_0,v_0)=\frac{\partial f_{\nu_0}}{\partial y}(u_0,v_0)=0,\quad 
\det(Hf_{\nu_0})_{(u_0,v_0)}\neq 0,
\end{equation}
where $(Hf_{\nu_0})_{(u_0,v_0)}$ is the Hessian matrix of $f_{\nu_0}$ at $(x,y)=(u_0,v_0)$.

The tangency $r$ \emph{unfolds generically} at $\nu = \nu_0$ if
$$\frac{\partial f}{\partial \nu}(\nu_0,u_0,v_0)\neq 0.$$
\end{definition}

\begin{remark}
It is easy to see that the property (1) does not depend on the coordinates used to 
set $l_\nu$ in the $z$-axis. 
Similarly, the properties (2) and (3) do not depend on the coordinates used to 
set $S_\nu$ in the $xy$-plane. 
\end{remark}

\begin{figure}[hbtp]
\centering
\scalebox{0.57}{\includegraphics[clip]{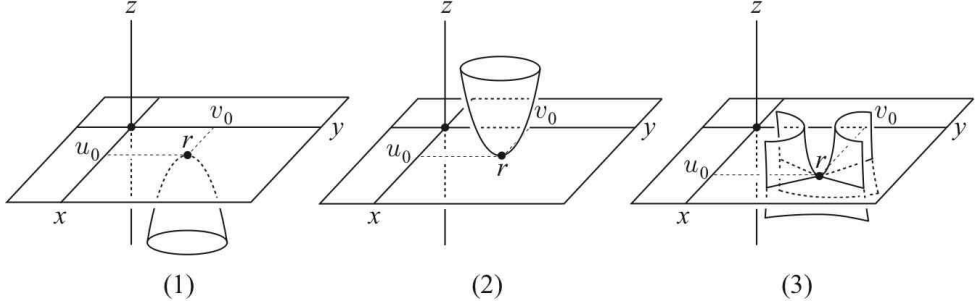}}
\caption{(1) $r$ is of elliptic type and ${\partial}^2 f_{\nu_0}(u_0,v_0)/{\partial} x^2<0$. 
(2) $r$ is of elliptic type and ${\partial}^2 f_{\nu_0}(u_0,v_0)/{\partial} x^2>0$.
(3) $r$ is of hyperbolic type.}
\label{fg_2}
\end{figure}

When $\det(Hf_{\nu_0})_{(u_0,v_0)}> 0$ (resp.\ $<0$) in Definition \ref{d_generic_unfolding}\,(3), we say that the 
tangency $r=(u_0,v_0,0)$ is of \emph{elliptic} (resp.\ \emph{hyperbolic}) 
type.
The Taylor expansion of $f_{\nu_0}$ around $(u_0,v_0)$ is 
\begin{equation}\label{form1}
\begin{split}
f_{\nu_0}(x,y)&=\frac{1}{2}\frac{{\partial}^2 f_{\nu_0}}{{\partial} x^2}(u_0,v_0) (x-u_0)^2+
\frac{{\partial}^2 f_{\nu_0}}{{\partial x}{\partial y} }(u_0,v_0)(x-u_0)(y-v_0)\\
 &\qquad+\frac{1}{2}\frac{{\partial}^2 f_{\nu_0}}{{\partial} y^2}(u_0,v_0)(y-v_0)^2+o\bigl((|x-u_0|+|y-v_0|)^2\bigr).
\end{split}
\end{equation}
From (\ref{form1}) together with the classification of quadratic surfaces in $\mathbb{R}^3$, we know that $Y_{\nu_0}$ has the 
form near $r=(u_0,v_0,0)$ as illustrated in Fig.\ \ref{fg_2}.

\subsection{Generic conditions}\label{generic_conditions}
Throughout the remainder of this paper, we suppose that 
$\varphi$ is a  $3$-dimensional $C^2$ diffeomorphism with saddle fixed points $p$ of $\mathrm{index}(p)=1$ and 
$q$ of  $\mathrm{index}(q)=2$, and such that $W^s(p)$ and $W^u(q)$ have a nondegenerate 
heterodimensional tangency $r$, $W^u(p)$ and $W^s(q)$ have a quasi-transverse intersection $s$.
The $\varphi$ is \emph{locally $C^2$ linearizable} in a neighborhood $U(q)$ of $q$ if 
there exists a $C^2$ linearizing coordinate $(x,y,z)$ on $U(q)$, that is,
\begin{equation}\label{linearizing}
q=(0,0,0),\quad \varphi(x,y,z)=({\alpha}x,{\beta}y,{\gamma}z)
\end{equation}
for any $(x, y, z)\in U(q)$ with $\varphi(x, y, z)\in U(q)$, where $\alpha, \beta$ and $\gamma$ are eigenvalues of $(d\varphi)_q$. 

One can take a local unstable manifold $W^u_\mathrm{loc}(q)$ so that it is an 
open disk in the plane $\{z=0\}$ centered at $(x,y)=(0,0)$.
We may assume that the both points $r,s$ are contained in $U(q)$ if necessary replacing $r$ (resp.\ $s$) 
by $\varphi^{-n}(r)$ (resp.\ $\varphi^n(s)$) with sufficiently large $n\in \mathbb{N}$.
We set 
$$r=(u_0,v_0,0)$$
with respect to the linearizing coordinate on $U(q)$.

We suppose moreover that $\{\varphi_{\mu, \nu}\}$  
is a two-parameter family in $\mathrm{Diff}^{2}(M)$ with $\varphi_{0,0}=\varphi$ and satisfying 
the conditions of Theorem \ref{thm_A}.
In particular, $\varphi_{\mu,\nu}$ is locally $C^2$ linearizable in a small neighborhood $U(q_{\mu,\nu})$ 
of $q_{\mu,\nu}$ in $M$ and hence $\varphi_{\mu,\nu}$ has the form as  (\ref{linearizing}) in $U(q_{\mu,\nu})$, where 
$\alpha,\beta,\gamma$ are $C^2$ functions on $\mu,\nu$, i.e.,\ $\alpha=\alpha_{\mu,\nu},\beta=\beta_{\mu,\nu},
\gamma=\gamma_{\mu,\nu}$.

We will put the following generic conditions (C1)-(C4) as the hypotheses in Theorem \ref{thm_A}.

\begin{enumerate}[({C}1)]
\item (Generic condition for $q$) 
The $\varphi$ is locally $C^2$ linearizable at $q$  given as in (\ref{linearizing}).
For simplicity,  we suppose that every eigenvalues  of $(d\varphi)_q$  is positive, that is, 
$$0<\gamma<1<\beta<\alpha.$$
\item (Generic unfolding property for $r$)
The nondegenerate heterodimensional tangency $r$ of 
$W^u(q)$ and $W^s(p)$ unfolds generically with respect to the $\mu$-parameter 
families $\{W^u(q_{\mu,0})\}$ and $\{W^s(p_{\mu,0})\}$.
\item (Generic unfolding property for $s$)
The quasi-transverse intersection $s$ of 
$W^s(q)$ and $W^u(p)$ unfolds generically with respect to the $\nu$-parameter 
families $\{W^s(q_{0,\nu})\}$ and $\{W^u(p_{0,\nu})\}$.
\item (Additional generic conditions)
The tangency $r$ is not on the $x$-axis $W^{uu}_\mathrm{loc}(q)$, that is,
\begin{equation}\label{generic conditions2}
v_0\neq 0.
\end{equation}

There exists a regular parametrization $l(t)= (x(t), y(t), z(t))$ $(t\in I)$ of a small curve in $W^u(p)\cap U(q)$ with respect to 
the linearizing coordinate $(x,y,z)$ on $U(q)$ with $s=l(0)$ and
\begin{equation}\label{generic conditions}
\frac{dx}{dt}(0)\neq 0,
\end{equation} 
where $I$ is an open interval centered at $0$.

There exists a $C^2$ function $f:O\rightarrow \mathbb{R}$ defined on an open disk $O$ in the $xy$-plane centered at 
$r$ such that $f(u_0,v_0)=0$,  $\{(x,y,f(x,y))\,;\,(x,y)\in O\}\subset W^s(p)\cap U(q)$ and 
\begin{equation}\label{generic conditions3}
\frac{{\partial}^2 f}{{\partial} x^2}(u_0,v_0)\neq 0.
\end{equation}
Note that the condition (\ref{generic conditions3}) is automatically satisfied when $r$ is of elliptic type.
\end{enumerate}

\section{Some lemmas about parametrization and curvatures}\label{some_lemma}
The goal of this section is to prove three lemmas needed for the proof of Theorem \ref{thm_A}.
These play important roles in Section \ref{Proof}.
\begin{itemize}

\item
Lemma \ref{new_parameters} presents a new parameter $(\hat{\mu}, \hat{\nu})$ such that, for 
any $\hat\mu$ near $0$, there exists a quasi-transverse intersection $s_{\hat{\mu}, 0}$ of $W^{s}(q_{\hat{\mu}, 0})$ and $W^{u}(p_{\hat{\mu}, 0})$ 
which unfolds generically at $\hat{\nu}=0$ with respect to the $\hat{\nu}$-parameter.  
After Lemma \ref{new_parameters}, we denote the new parameter $(\hat{\mu}, \hat{\nu})$ again by $(\mu, \nu)$ 
for simplicity.

\item
In Lemma \ref{inclination}, we show that, for any $\mu_0$ near $0$, there exists a regular curve $l_{m}$ in $W^{u}(p_{{\mu}_{0}, 0}) $ containing the quasi-transverse
intersection ${\varphi}^{m}_{{\mu}_{0}, 0}(s_{{\mu}_{0}, 0})$ and arbitrarily $C^2$ close to $W^{uu}_{\mathrm{loc}}(q_{{\mu}_{0}, 0})$.
In particular, this implies that the curvature of $l_{m}$ can be taken arbitrarily close to $0$ with respect to the linearizing coordinate (\ref{linearizing}) on $U(q_{\mu_0,0})$.

\item
Lemma \ref{quadratic proof} gives a connection between the curvature and quadratic tangencies.
In fact, we show that a tangency $\tau$ of a regular curve $l$ and a regular surface $S$ in $\mathbb{R}^{3}$ 
is quadratic if the curvature of $l$ at $\tau$ is different from the normal curvature of $S$ at $\tau$ 
along the direction tangent to $l$.
\end{itemize}

\bigskip

\begin{figure}[hbtp]
\centering
\scalebox{0.57}{\includegraphics[clip]{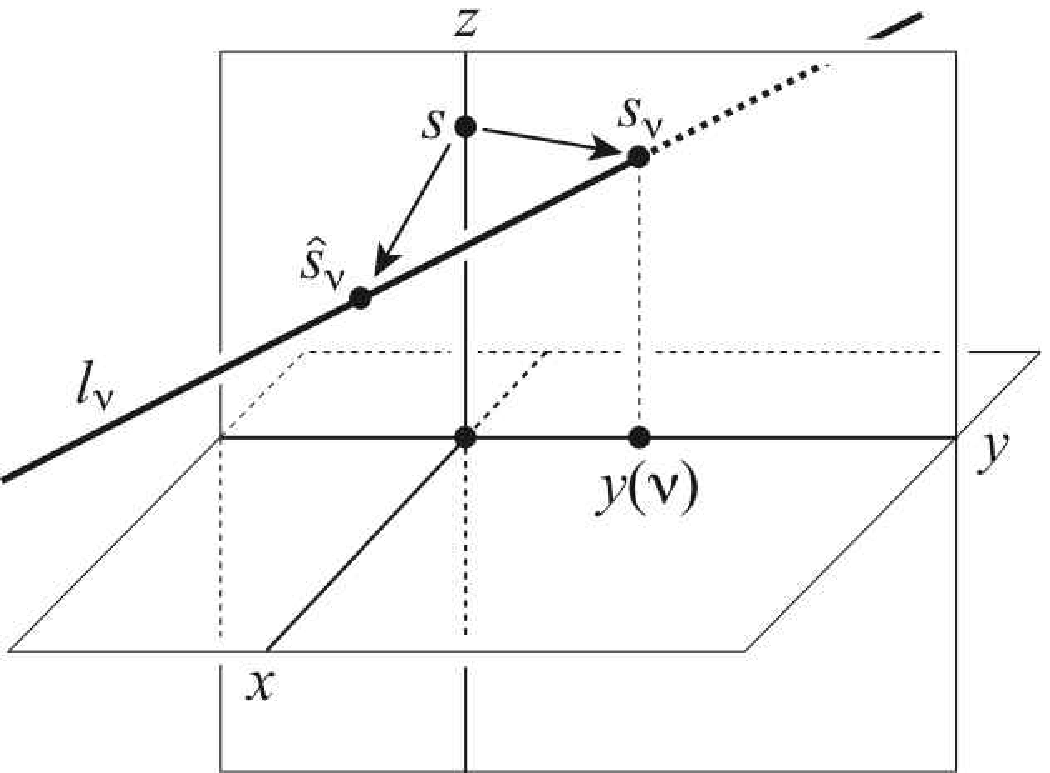}}
\caption{}
\label{fg_3}
\end{figure}

For any $(\mu,\nu)$ near $(0,0)$, we may assume that $U(q_{\mu,\nu})$ is equal to 
$$D(\delta):=(-\delta,\delta)^3$$
with respect to the linearizing coordinate given in Subsection \ref{generic_conditions} for some 
constant $\delta>0$.
Since $s$ is a quasi-transverse intersection which unfolds generically with respect to the $\nu$-parameter families 
$\{W^s(q_{0,\nu})\}$ and $\{W^u(p_{0,\nu})\}$ by the condition (C3), there exists a $C^2$ continuation $\hat s_\nu\in 
W^u(p_{0,\nu})\cap D(\delta)$ with $\hat s_{0}=s$ and such that $\hat s_{\nu}$ satisfies the 
conditions same as those for $s_\nu$ in Definition \ref{d_generic_unfolding}\,(1).
By (\ref{generic conditions}), for any $\nu$ near $0$, the component $l_\nu$ of $W^u(p_{0,\nu})\cap D(\delta)$ containing 
$\hat s_\nu$ meets 
transversely the $yz$-plane at a point $s_\nu$ which defines a $C^2$ continuations $\{s_\nu\}$ with $s_0=s$, see Fig.\ \ref{fg_3}.
Note that $d \hat s_\nu/d\nu(0)=d s_\nu/d\nu(0)+\boldsymbol{w}$ for some $\boldsymbol{w}\in T_s(l_0)=T_s(W^u(p))$, 
where $d \hat s_\nu/d\nu(0)$ denotes $d \hat s_\nu/d\nu|_{\nu=0}$.
Let $y(\nu)$ be the $y$-coordinate of $s_\nu$.
If $d y/d\nu(0)=0$, then $d s_\nu/d\nu(0)$ would be tangent to the $z$-axis $W_{\mathrm{loc}}^s(q)$ at $s$ 
and hence $d \hat s_\nu/d\nu(0)\in T_s(W^s(q))\oplus T_s(W^u(p))$.
This contradicts (\ref{quasi_transverse}).
Thus, we have 
\begin{equation}\label{d_s}
\frac{d y}{d\nu}(0)\neq 0.
\end{equation}
For any $(\mu,\nu)$ near $(0,0)$, let $s_{\mu,\nu}$ be a transverse intersection point of $W^u(p_{\nu,\nu})\cap D(\delta)$ 
with the $yz$-plane such that $\{s_{\mu,\nu}\}$ is a $C^2$ continuation with $s_{0,\nu}=s_\nu$.

\begin{lemma}\label{new_parameters}
There exists a constant  $\rho>0$ and  a $C^2$ function
$ \tilde{\nu} :(-\rho, \rho)\rightarrow \mathbb{R}$ 
such that, for any $\mu \in (-\rho, \rho)$,  
$s_{\mu,\tilde{\nu}(\mu)}$ is a quasi-transverse intersection of $W^s(q_{\mu, \tilde{\nu}(\mu)})$ and $W^u(p_{\mu, \tilde{\nu}(\mu)})$ 
which unfolds generically  with respect to the $\nu$-parameter families  $\{W^s(q_{\mu(\mathrm{fixed}), \nu})\}$ and 
$\{W^u(p_{\mu(\mathrm{fixed}), \nu})\}$. 
\end{lemma}

\begin{proof}
Let $y(\mu,\nu)$ be the $y$-coordinate of $s_{\mu,\nu}$.
By (\ref{d_s}), $\partial y/\partial \nu(0,0)\neq 0$.
Hence, 
by the Implicit Function Theorem, 
there exists 
a $C^2$ function $\tilde{\nu}: (-\rho, \rho)\rightarrow \mathbb{R}$ for some  $\rho>0$ 
such that 
$$\tilde{\nu}(0)=0,\  
y(\mu, \tilde{\nu}(\mu))=0,\ \frac{\partial y}{\partial \nu}(\mu, \tilde\nu(\mu))\neq 0
$$
for any $\mu\in  (-\rho, \rho)$.
This implies that $s_{\mu,\tilde\nu(\mu)}$ 
is a quasi-transverse intersection
unfolding generically at $\nu=\tilde\nu(\mu)$  with respect to the $\nu$-parameter families 
$\{W^s(q_{\mu, \nu})\}$ and $\{W^u(p_{\mu, \nu})\}$. 
\end{proof}

\subsection*{A new parametrization}
Consider the coordinate $(\hat \mu,\hat \nu)$ on the parameter space defined by $\hat\mu=\mu, \hat\nu=\nu-
\tilde\nu(\mu)$.
For simplicity, we denote the new coordinate again by $(\mu,\nu)$.
Thus, there exists a continuation $\{s_{\mu,0}\}_{\mu\in (-\rho,\rho)}$ of quasi-transverse intersections 
of $W^s(q_{\mu,0})$ and $W^u(p_{\mu, 0})$ such that each $s_{\mu,0}$ unfolds generically at $\nu=0$ with 
respect to the $\nu$-parameter families $\{W^s(q_{\mu,\nu})\}$ and 
$\{W^u(p_{\mu, \nu})\}$.

\bigskip

Fix $\mu_0$ with sufficiently small $|\mu_0|$ arbitrarily.
By the properties (\ref{linearizing}) and (\ref{generic conditions}), there exists $m_0\in \mathbb{N}$ such that, 
for any $m\geq m_0$, one can parameterize the component 
$l_m$ of $W^u(p_{\mu_0,0})\cap D(\delta)$ containing $\varphi^m_{\mu_0,0}(s_{\mu_0,0})$ so that 
$l_m(0)=\varphi^m_{\mu_0,0}(s_{\mu_0,0})$ and 
$$l_m(t)=(t,y_m(t),z_m(t))\quad (t\in (-\delta,\delta)).$$

\begin{lemma} \label{inclination}
The sequence $\{l_{m}\}$ $C^2$ converges uniformly to $W^{uu}_{\mathrm{loc}}(q_{\mu_0,0})$ as $m\rightarrow \infty$.
In particular, for any $\varepsilon >0$, there exists $\hat m_{0}\geq m_0$ such that  the curvature 
at any point of $l_{m}$ is less than $\varepsilon$ with respect to the standard Euclidean metric on 
$U(q_{\mu_0,0})=D(\delta)$ if $m\geq \hat m_0$.
\end{lemma}
\begin{proof}
By (\ref{linearizing}), for any $m\geq m_0$, 
$$
l_{m}(t)=(t, {\beta}^n y_{m_0}(\alpha^{-n}t), {\gamma}^n z_{m_0} (\alpha^{-n}t)),
$$
where $n=m-m_0$, $\alpha=\alpha_{\mu_0,0}$, $\beta=\beta_{\mu_0,0}$, $\gamma=\gamma_{\mu_0,0}$.
Thus we have
\begin{equation}\label{limits}
\begin{split}
\frac{dl_{m}}{dt}(t)&=\left(1,\ \frac{\beta^n}{\alpha^n}\frac{dy_{m_0}}{dt}(\alpha^{-n}t),\ 
\frac{\gamma^n}{\alpha^n}\frac{dz_{m_0}}{dt}(\alpha^{-n}t)\right)\xrightarrow{\mathrm{uniformly}} (1,0,0),\\
\frac{d^2l_{m}}{dt^2}(t)&=\left(0,\ \frac{\beta^n}{\alpha^{2n}}\frac{d^2y_{m_0}}{dt^2}(\alpha^{-n}t),\ 
\frac{\gamma^n}{\alpha^{2n}}\frac{d^2z_{m_0}}{dt^2}(\alpha^{-n}t)\right)\xrightarrow{\mathrm{uniformly}} (0,0,0)
\end{split}
\end{equation}
as $m\rightarrow \infty$.
Since $\{l_m(0)\}_{m=m_0}^\infty$ converges to $q_{\mu_0,0}=(0,0,0)$, it follows from (\ref{limits}) that $\{l_m\}$ 
$C^2$ converges uniformly to the $x$-axis in $D(\delta)$. 
\end{proof}  

\begin{lemma}\label{quadratic proof}
Let $S$ be a regular surface in the Euclidean $3$-space $\mathbb{R}^3$ and $l$ a regular curve tangent to $S$ at $\tau$.
Suppose that the curvature $\kappa_l(\tau)$ of $l$ at $\tau$ is less than the absolute value of 
the normal curvature $\kappa_S(\tau,\boldsymbol{w})$ of $S$ at $\tau$ along a non-zero vector 
$\boldsymbol{w}$ tangent to $l$.
Then tangency of $S$ and $l$ at $\tau$ is quadratic.
\end{lemma}
\begin{proof}
By changing the coordinate $(x,y,z)$ on $\mathbb{R}^3$ by an isometry, we may assume that $\tau=(0,0,0)$, the tangent 
space of $S$ at $\tau$ is the $xy$-plane and $\boldsymbol{w}/\Vert \boldsymbol{w}\Vert=(1,0,0)$.
Then one can suppose that $S$ (resp.\ $l$) is parameterized as $(x,y,\psi(x,y))$ (resp.\ 
$(x,f_1(x),f_2(x))$ in a small neighborhood of $(0,0,0)$ in $\mathbb{R}^3$.
Since the graph of $z=\psi(x,0)$ is the cross section of $S$ along the $xz$-plane, 
$$|\kappa_S(\tau,\boldsymbol{w})|=\frac{|g''(0)|}{(g'(0)^2+1)^{3/2}}=|g''(0)|,$$
where $g(x)=\psi(x,0)$.
Since the graph of $z=f_2(x)$ coincides with the orthogonal projection $\overline l$ of $l$ into the $xz$-plane,
$$\kappa_l(\tau)\geq \kappa_{\overline l}(\tau)=\frac{|f_2''(0)|}{(f_2'(0)^2+1)^{3/2}}=|f_2''(0)|.$$
It follows from our assumption $|\kappa_S(\tau,\boldsymbol{w})|>\kappa_l(\tau)$ that 
$|g''(0)|>|f_2''(0)|$.
This shows that the tangency at $\tau$ is quadratic.
\end{proof}

\section{Proof of Theorem \ref{thm_A}}\label{Proof}
In this section, we give the proof of Theorem \ref{thm_A}.
\begin{itemize}
\item In Subsection \ref{Existence},
we show that, for any $\mu_0$ in either $(-\varepsilon,0)$ or $(0,\varepsilon)$ and any sufficiently large $m\in \mathbb{N}$, 
there exists $\nu_m$ with $\lim_{m\rightarrow \infty}\nu_m=0$ such that $W^{u}(p_{{\mu}_{0}, \nu_m})$ and $W^{s}(p_{{\mu}_{0}, \nu_m})$ have a quadratic tangency $\tau_m$ (Assertion \ref{elliptic} and Assertion \ref{hyperbolic}).
Here the sign of $\mu_0$ is chosen so that $\mu_0\cdot b_{\mu_0,0}<0$ (resp.\ $\mu_0\cdot b_{\mu_0,0}>0)$ if the tangency $r$ is 
of elliptic (resp.\ hyperbolic) type, where $b_{\mu_0,0}$ is the coefficient of $(x-\nu_{\mu_0,0})^2$-term of the Taylor expansion (\ref{form4}).
See Fig.\ \ref{fg_4} and Fig.\ \ref{fg_7}.
As is suggested in Fig.\ \ref{fg_5}, the existence of the homoclinic tangency $\tau_m$ is proved by using the Intermediate Value Theorem.
By Lemma \ref{inclination}, the curvature of $W^u(p_{\mu_0,\nu_m})$ at $\tau_m$ converges to zero as $m\rightarrow \infty$.
On the other hand, we will show that the normal curvature of $W^s(p_{\mu_0,\nu_m})$ at $\tau_m$ along the direction tangent to $l_m$ is bounded away from zero.
Hence, by Lemma \ref{quadratic proof}, the tangency $\tau_m$ is quadratic.
 
\item 
In Subsection \ref{Unfolding generically}, we show that the quadratic homoclinic tangency ${\tau}_{m}$ unfolds generically at $\nu={\nu}_{m}$ with respect to the $\nu$-parameter families $\{W^{s}(p_{{\mu}_{0} ,\nu})\}$ and $\{W^{u}(p_{{\mu}_{0}, \nu})\}$ by representing a 
neighborhood of $\tau_m$ in $W^s(p_{\mu_0,\nu})$ as the graph of a function  of $(x,z)$.
\end{itemize}

\subsection{Existence of quadratic homoclinic tangencies}\label{Existence}
Let $\{\varphi_{\mu, \nu}\}$ be the family given in Subsection \ref{generic_conditions}.
In particular, $r=(u_0,v_0,0)$ is a nondegenerate heterodimensional tangency of $W^u(q)$ and 
$W^s(p)$ which unfolds generically with respect to the $\mu$-parameter families 
$\{W^u(q_{\mu,0})\}$ and $\{W^s(p_{\mu,0})\}$.
By our settings in Sections \ref{Generic_conditions} and \ref{some_lemma}, there exist $C^2$ functions 
$f_{\mu,\nu}:O\subset \mathbb{R}^2\rightarrow \mathbb{R}$ $C^2$ 
depending on $(\mu,\nu)$ with $f_{0,0}=f$ and 
$$\Sigma(\mu,\nu):=\{(x,y,f_{\mu,\nu}(x,y))\,;\,(x,y)\in O\}\subset W^s(p_{\mu,\nu})\cap D(\delta)$$
for any $(\mu,\nu)$ near $(0,0)$.
Since $\det(Hf)_{(u_0,v_0)}\neq 0$, there exists a uniquely determined $C^2$ continuation $(u_{\mu,\nu},v_{\mu,\nu})$ 
with $(u_{0,0}, v_{0,0})=(u_0,v_0)$ and 
$$\frac{\partial f_{\mu,\nu}}{\partial x}(u_{\mu,\nu},v_{\mu,\nu})=\frac{\partial f_{\mu,\nu}}{\partial y}(u_{\mu,\nu},v_{\mu,\nu})=0.$$

\begin{proposition}\label{prop1}
For a sufficiently small $\varepsilon >0$ and any $\mu$ in either $(0,\varepsilon)$ or $(-\varepsilon,0)$, 
there exists $\nu$ arbitrarily close to $0$ such that $\varphi_{\mu, \nu}$ has a quadratic homoclinic tangency 
associated to $p_{\mu, \nu}$.  
\end{proposition}

By the condition (\ref{generic conditions2}), $r$ is not in the $x$-axis.
One can take the linearizing coordinate on $D(\delta)$ so that $s$ (resp.\ $r$) is in the upper half space $\{z>0\}$ 
(resp.\ $\{x>0\}$).
The Taylor expansion of $f_{\mu,\nu}$ around $(u_{\mu,\nu},v_{\mu,\nu})$ has the form: 
\begin{equation}\label{form4}
\begin{split}
f_{\mu,\nu}(x,y)=a_{\mu,\nu}&+\frac{1}{2}b_{\mu,\nu}(x-u_{\mu,\nu})^2+c_{\mu,\nu}(x-u_{\mu,\nu})(y-v_{\mu,\nu})\\
&+\frac{1}{2}d_{\mu,\nu}(y-v_{\mu,\nu})^2+o\bigl((|x-u_{\mu,\nu}|+|y-v_{\mu,\nu}|)^2\bigr),
\end{split}
\end{equation}
where $a_{0,0}=0$ and 
$$b_{\mu,\nu}=\frac{{\partial}^2 f_{\mu,\nu}}{{\partial} x^2}(u_{\mu,\nu},v_{\mu,\nu}),\ 
c_{\mu,\nu}=\frac{{\partial}^2f_{\mu,\nu}}{{\partial x}{\partial y}}(u_{\mu,\nu},v_{\mu,\nu}),\ 
d_{\mu,\nu}= \frac{{\partial}^2 f_{\mu,\nu}}{{\partial} y^2}(u_{\mu,\nu},v_{\mu,\nu}).$$
Since the tangency $r$ unfolds generically with respect to $\varphi=\varphi_{0,0}$ by (C1), 
\begin{equation}\label{e_alpha}
\eta_0=\frac{\partial a_{\mu,\nu}}{\partial \mu}\Big|_{(\mu,\nu)=(0,0)}\neq 0.
\end{equation}
If necessary replacing $\mu$ by $-\mu$, we may assume that $\eta_0>0$.
By the condition (\ref{generic conditions3}), $b_{0,0}\neq 0$ and hence $b_{\mu,\nu}\neq 0$ for any $(\mu,\nu)$ near 
$(0,0)$. 

Proposition \ref{prop1} is divided to the following two assertions.

\begin{assertion}[Elliptic case]\label{elliptic}
If $r$ is of elliptic type, then Proposition \ref{prop1} holds.
\end{assertion}
\begin{proof}
First we consider the case of $b_{\mu,\nu}<0$ for any $(\mu,\nu)$ near $(0,0)$.
By (\ref{e_alpha}), for any sufficiently small $\mu_0>0$, the intersection 
$C_{\mu_0}=\Sigma(\mu_0,0)\cap \{z=0\}$ is a circle disjoint from the $x$-axis.
For a sufficiently small $h_0>0$, $A=\Sigma(\mu_0,0)\cap \{0\leq z\leq h_0\}$ is 
an annulus in $D(\delta)$, see Fig.\ \ref{fg_4}\,(1).
\begin{figure}[hbtp]
\centering
\scalebox{0.6}{\includegraphics[clip]{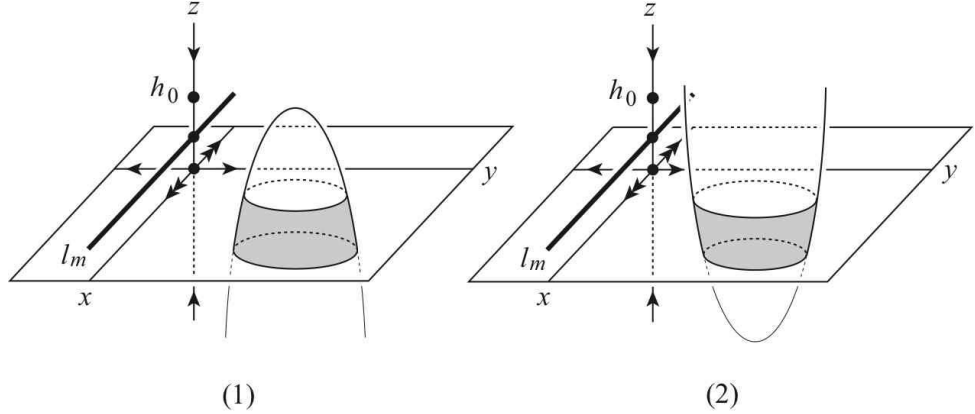}}
\caption{(1) The case of $\mu_0>0,\nu=0, b_{\mu_0,0}<0$. (2) The case of $\mu_0<0,\nu=0,b_{\mu_0,0}>0$.
Each shaded region represents $A$.}
\label{fg_4}
\end{figure}
Replacing $m_0$ by an integer greater than $m_0$ if necessary, we may assume that $z_m(0)<h_0/2$ for any $m\geq m_0$.
By Lemma \ref{inclination}, the curve $l_m\subset W^u(p_{\mu_0,0})\cap D(\delta)$ given in Section \ref{some_lemma} 
is sufficiently $C^2$ close to the $x$-axis.
Thus one can suppose that $\pi_y(l_{m})\cap \pi_y(A)=\emptyset$, where $\pi_y:D(\delta)\longrightarrow \mathbb{R}$ 
is the orthogonal projection defined by $\pi_y(x,y,z)=y$.
For any sufficiently small $\nu$, let $l_{m,\nu}$ be the component of $W^u(p_{\mu_0,\nu})\cap D(\delta)$ such that 
$\{l_{m,\nu}\}$ 
is an $\nu$-continuation with $l_{m,0}=l_{m}$, and set $A_\nu=\Sigma({\mu_0,\nu})\cap \{0\leq z\leq h_0\}$.
Moreover, one can suppose that $l_{m,\nu}$ is parameterized as
$$l_{m,\nu}(t)=(t,y_m(\nu,t),z_m(\nu,t))\quad (t\in (-\delta,\delta)).$$
By the condition (C3), one can take $\bar\nu\neq 0$ with arbitrarily small $|\bar\nu|$ such that 
$$0<\pi_y(l_{m_0,\bar\nu}(0))\leq \sup\{\pi_y(l_{m_0,\bar\nu})\}<\min\{\pi_y(A_{\bar\nu})\}.$$
We may assume that $\bar\nu>0$ if necessary replacing $\nu$ by $-\nu$.
For any integer $m$ sufficiently greater than $m_0$, there exists $0<\bar \nu_m<\bar \nu$ such 
the continuation $\{l_{m,\nu}\}_{0 \leq \nu \leq \bar \nu_{m}}$ is well defined and 
$$\max\{\pi_y(A_{\bar\nu_m})\}<\inf\{\pi_y(l_{m,\bar\nu_m})\}$$
holds, see Fig.\ \ref{fg_5}\,(1).
By the Intermediate Value Theorem, there exists $0<\nu_m<\bar\nu_m$ such that 
$l_{m,\nu_m}$ and $A_{\nu_m}$ have a tangency $\tau_m$, see Fig.\ \ref{fg_5}\,(2).
\begin{figure}[hbtp]
\centering
\scalebox{0.62}{\includegraphics[clip]{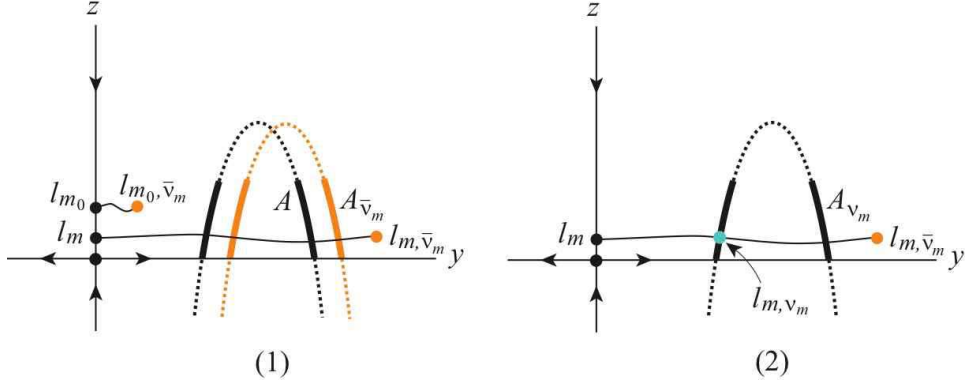}}
\caption{The cross sections.}
\label{fg_5}
\end{figure}
Since $l_{m,\nu_m}\subset W^u(p_{\mu_0,\nu_m})$ and $A_{\nu_m}\subset W^s(p_{\mu_0,\nu_m})$, 
$\tau_m$ is a homoclinic tangency associated to $p_{\mu_0,\nu_m}$.

When $b_{\mu,\nu}>0$ for any $(\mu,\nu)$ near $(0,0)$, one can prove the existence of a homoclinic tangency 
$\tau_m$ near $r$ associated to $p_{\mu_0,\nu_m}$ 
by arguments quite similar to those as above for any $\mu_0$ with $\mu_0<0$.

It remains to show that the tangency $\tau_m$ is quadratic.
Since $\Sigma(\mu_0,\nu_m)$ is of elliptic type and $\lim_{m\rightarrow \infty}\nu_m=0$, any normal curvature of 
$\Sigma(\mu_0, \nu_m)$ at $\tau_m$ is greater than some positive constant $\kappa_0$ independent of $m$.
On the other hand, by an argument quite similar to that in Lemma \ref{inclination}, for any $m$ sufficiently greater than 
$m_0$, the curvature of $l_{m,\nu_m}$ at $\tau_m$ is less than $\kappa_0$.
Thus, by Lemma \ref{quadratic proof}, $\tau_m$ is a quadratic tangency.
\end{proof}

\begin{assertion}[Hyperbolic case]\label{hyperbolic}
When  $r$ is a tangency of hyperbolic type, Proposition \ref{prop1} holds.
\end{assertion}

\begin{proof}
Since $r$ is of hyperbolic type, $\Sigma(0,0)\cap \{z=0\}$ consists of two almost 
straight curves $\alpha_1,\alpha_2$ meeting transversely at $r$, see Fig.\ \ref{fg_6}\,(1).
\begin{figure}[hbtp]
\centering
\scalebox{0.6}{\includegraphics[clip]{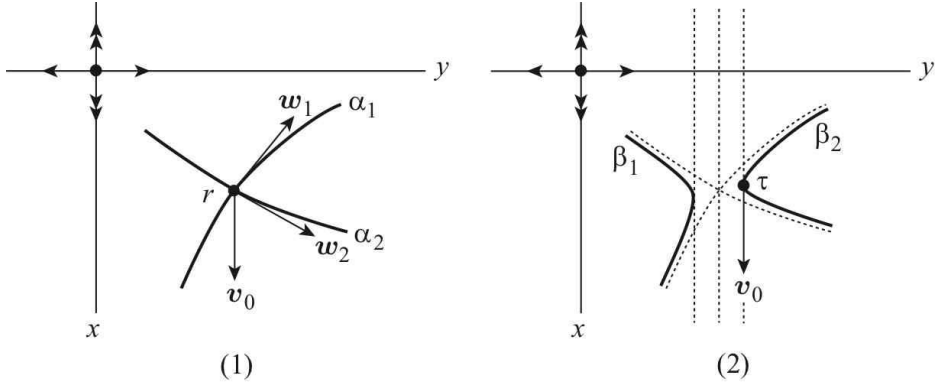}}
\caption{(1) The situation when $(\mu,\nu)=(0,0)$.
(2) The situation when $(\mu,\nu)=(\mu_0,0)$.}
\label{fg_6}
\end{figure}
If necessary by reducing the domain $O$ of $f_{\mu,\nu}$ containing $(u_{\mu,\nu},v_{\mu,\nu})$, we may assume that 
$\Sigma(\mu,\nu)\cap \{z=0\}$ is disjoint from the $x$-axis for any $(\mu,\nu)$ near $(0,0)$.
If $\boldsymbol{w}_i=(\xi_i,\eta_i,0)$ $(i=1,2)$ is a unit vector tangent to $\alpha_i$ at $r$, then 
$b_{0,0}\xi_i^2+2c_{0,0}\xi_i\eta_i+d_{0,0}\eta_i^2=0$.
This implies that the normal curvature $\kappa_{\Sigma(0,0)}(r,\boldsymbol{w}_i)$ of $\Sigma(0,0)$ at $r$ 
along $\boldsymbol{w}_i$ is zero.
Since $b_{0,0}\neq 0$ by (\ref{generic conditions3}), both $\boldsymbol{w}_1$, $\boldsymbol{w}_2$ are 
not parallel to the unit tangent vector $\boldsymbol{v}_0=(1,0,0)$.
Thus we have $\kappa_{\Sigma(0,0)}(r,\boldsymbol{v}_0)\neq 0$.
When $b_{0,0}<0$ (resp.\ $b_{0,0}>0$), for any sufficiently small $\mu_0$ with $\mu_0<0$ (resp.\ $\mu_0>0$), 
$\Sigma(\mu_0,0)\cap \{z=0\}$ consists of two $C^2$ curves $\beta_1,\beta_2$ separated by a line in the $xy$-plane 
parallel to $x$-axis, see Fig.\ \ref{fg_6}\,(2), and 
\begin{equation}\label{e_kappa}
\kappa_0:=|\kappa_{\Sigma(\mu_0,0)}(\tau,\boldsymbol{v}_0)|>0,
\end{equation}
where $\tau$ is a point of $\beta_1\cup \beta_2$ the tangent line at which is parallel to the 
$x$-axis.
One can take $\bar\nu>0$ and $h_0>0$ so that $A_{\nu}=\Sigma(\mu_0,\nu)\cap \{0\leq z\leq h_0\}$ is 
a disjoint union of two curvilinear rectangles for any $\nu$ with $0\leq \nu\leq \bar\nu$, see Fig.\ \ref{fg_7}.
\begin{figure}[hbtp]
\centering
\scalebox{0.6}{\includegraphics[clip]{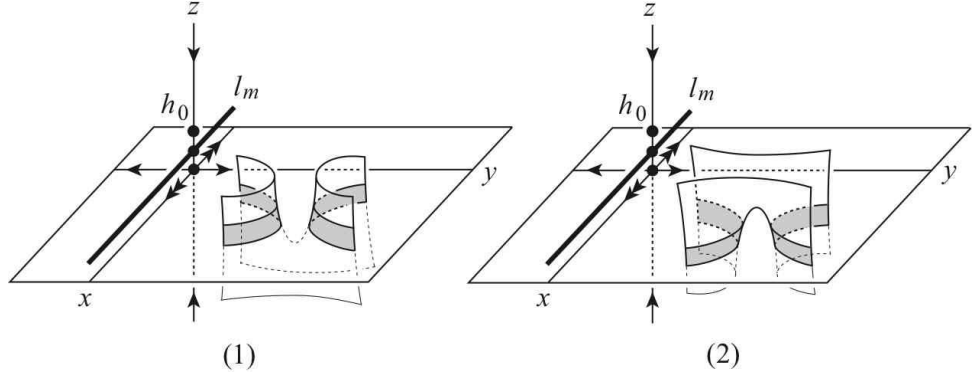}}
\caption{(1) The case of $\mu_0<0,\nu=0,b_{\mu_0,0}<0$. (2) The case of $\mu_0>0,\nu=0,b_{\mu_0,0}>0$.
Each pair of the shaded regions represents $A$.}
\label{fg_7}
\end{figure}
Moreover, by (\ref{e_kappa}), the $\bar\nu>0$ can be chosen so that 
$|\kappa_{\Sigma(\mu_0,\nu)}(\tilde\tau,\boldsymbol{w})|>\kappa_0/2$ for any point $\tilde\tau\in A_{\mu_0,\nu}$ sufficiently 
near $\tau$ and any unit vector $\boldsymbol{w}\in T_{\tilde\tau}(A_{\mu_0,\nu})$ sufficiently near $\boldsymbol{v}_0$.
As in the proof of Assertion \ref{elliptic}, for any integer $m$ sufficiently greater than $m_0$, there exists $\nu_m$ with 
$0<\nu_m<\bar\nu$ and
such that $l_{m,\nu_m}\subset W^s(p_{\mu_0,\nu_m})$ and $A_{\nu_m}\subset W^u(p_{\mu_0,\nu_m})$ 
have a quadratic tangency $\tau_m$, see Fig.\ \ref{fg_8}.
\begin{figure}[hbtp]
\centering
\scalebox{0.6}{\includegraphics[clip]{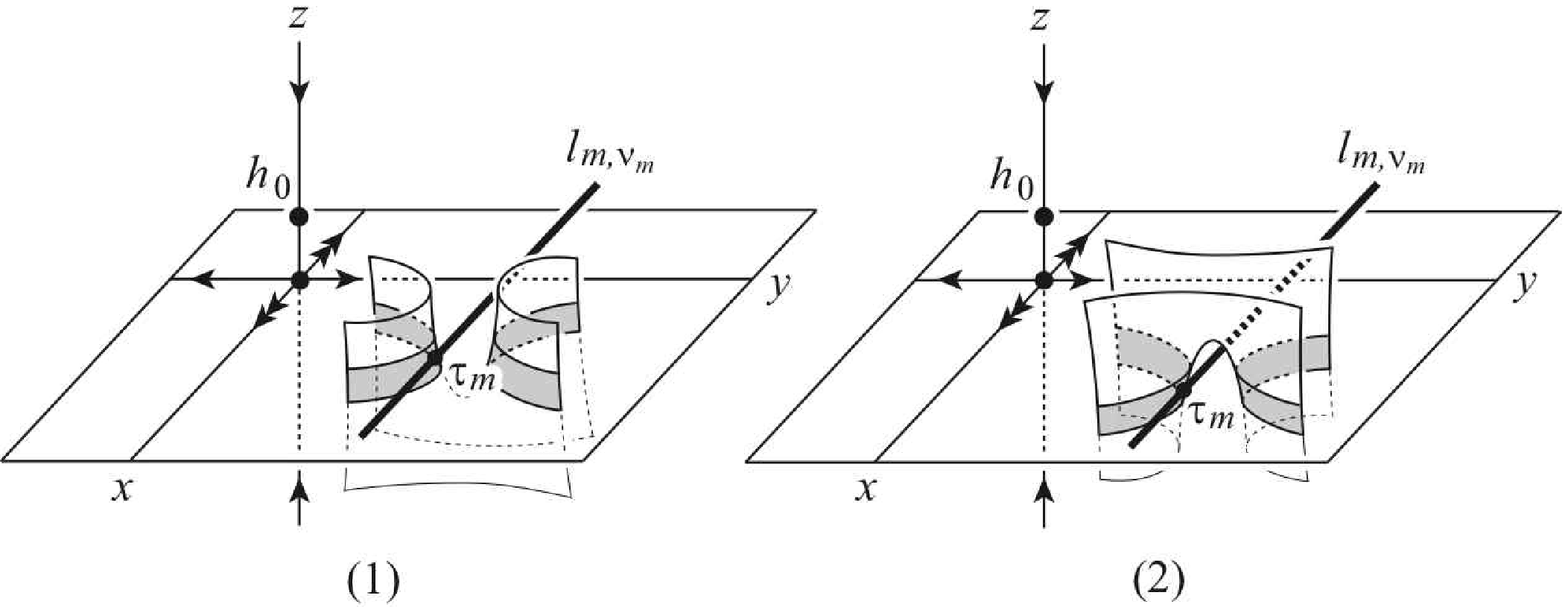}}
\caption{}
\label{fg_8}
\end{figure}
\end{proof}

\subsection{Generic unfolding of the tangency}\label{Unfolding generically}
For short, set $p_{\mu_0,\nu}=p_\nu$, $f_{\mu_0,\nu}(x,y)=f_\nu(x,y)$ 
and $(u_{\mu_0,\nu},v_{\mu_0,\nu})=(u_{\nu},v_{\nu})$.

Let $\tau_m=(\hat x_m,\hat y_m,f_{\nu_m}(\hat x_m,\hat y_m))$ be the homoclinic tangency of 
$W^{u}(p_{\nu_m})$ and 
$W^{s}(p_{\nu_m})$ given in Proposition \ref{prop1}.
From (\ref{form4}), we have
\begin{align*}
\frac{\partial f_{\nu_m}}{\partial y}(x,y)&=c_m(x-u_{\nu_m})+d_m(y-v_{\nu_m})+o_1,\\
\frac{\partial f_{\nu_m}}{\partial x}(x,y)&=b_m(x-u_{\nu_m})+c_m(y-v_{\nu_m})+o_1,
\end{align*}
where $b_m=b_{\mu_0,\nu_m}$, $c_m=c_{\mu_0,\nu_m}$, $d_m=d_{\mu_0,\nu_m}$ and $o_1=o(|x-u_{\nu_m}|+|y-v_{\nu_m}|)$.
Thus $b_m\partial f_{\nu_m}(x,y)/{\partial y}-c_m\partial f_{\nu_m}(x,y)/{\partial x}=(b_md_m-c_m^2)(y-v_{\nu_m})+o_1$.
On the other hand, since there exists a unit vector tangent to $\Sigma(\mu_0,\nu_m)$ at $\tau_m$ 
converges to $(1,0,0)$ as $m\rightarrow \infty$, $\lim_{m\rightarrow \infty}\partial f_{\nu_m}
(\hat x_m,\hat y_m)/{\partial x}=0$.
Since $\lim_{m\rightarrow \infty}b_m=b_{\mu_0,0}\neq 0$ and $\lim_{m\rightarrow \infty}b_md_m-c_m^2=\det(Hf_{\mu_0,0})
(u_{\mu_0,0},v_{\mu_0,0})\neq 0$,  
$$\frac{\partial f_{\nu_m}}{\partial y}(\hat x_m,\hat y_m)=\frac{c_m}{b_m}\frac{\partial f_{\nu_m}}{\partial x}
(\hat x_m,\hat y_m)+\frac{b_md_m-c_m^2}{b_m}\,(\hat y_m-v_{\nu_m})+o_1\neq 0$$
for all sufficiently large $m$.
By the Implicit Function Theorem, there exists a $C^2$ function $y=g_\nu(x,z)=g(\nu,x,z)$ defined in a small 
neighborhood of $(\nu_m,\hat x_m,f_{\nu_m}(\hat x_m,\hat y_m))$ in the $(\nu,x,z)$-space with
$$(x,y,f_\nu(x,y))=(x,g_\nu(x,z),z).$$

\begin{proposition}\label{prop2}
For all sufficiently large $m$, 
the quadratic homoclinic tangency $\tau_m$ of $W^{s}(p_{\nu_m})$ and $W^{u}(p_{\nu_m})$
unfolds generically at $\nu=\nu_m$ with respect to the $\nu$-parameter families $\{W^{s}(p_{\nu})\}$ and 
$\{W^{u}(p_{\nu})\}$.  
\end{proposition}
\begin{proof}
Recall that $l_{m,\nu}$ has the parametrization $l_{m,\nu}(t)=(t,y_{m}(\nu,t),z_m(\nu,t))$  
with $l_{m,\nu_m}(\hat x_m)=\tau_m$.
By Definition \ref{d_generic_unfolding}\,(2), it suffices to show that
\begin{equation}\label{eqn_1}
\frac{\partial y_{m}}{\partial \nu}(\nu_m,\hat x_m)\neq 
\frac{\partial g}{\partial \nu}(\nu_m,\hat x_m,z_m(\nu_m,\hat x_m))
\end{equation}
for all sufficiently large $m$.
Note that
$$\lim_{m\rightarrow \infty}\frac{\partial g}{\partial \nu}(\nu_m,\hat x_m,z_m(\nu_m,\hat x_m))=\frac{\partial g}{\partial \nu}
(0,\hat x_\infty,0),$$
where $\hat x_\infty$ is the $x$-coordinate of a  
point $\tau$ in $\Sigma(\mu_0,0)\cap \{z=0\}$ the tangent line in $xy$-plane at which is parallel to $(1,0,0)$, 
see Fig.\ \ref{fg_6}\,(2) in the case that $r$ is of hyperbolic type.
If we set $\tilde x_{m,\nu}=\alpha_{\nu}^{-n}\hat x_m$, then 
$\varphi_\nu^{n}(l_{m_0,\nu}(\tilde x_{m,\nu}))=l_{m,\nu}(\hat x_m)$, where $n={m-m_0}$ and 
$\alpha_\nu=\alpha_{\mu_0,\nu}$.
As was seen in the proof of Lemma \ref{new_parameters},
\begin{equation}\label{eqn_2}
\lim_{m\rightarrow \infty}\frac{\partial y_{m_0}}{\partial \nu}(\nu_m,\tilde x_{m,\nu_m})=
\frac{\partial y_{m_0}}{\partial \nu}(0,0)\neq 0.
\end{equation}
We denote the $\nu$-function $y_{m_0}(\nu,\tilde x_{m,\nu})$ by $h_m(\nu)$.
Since $\lim_{m\rightarrow \infty}d\tilde x_{m,\nu}/d\nu=0$,  
 it follows from (\ref{eqn_2}) that 
\begin{equation}\label{e_c}
\begin{split}
\left|\frac{d h_m}{d \nu}(\nu_m)\right|&=
\left|\frac{\partial y_{m_0}}{\partial \nu}(\nu_m,\tilde x_{m,\nu_m})+
\frac{\partial y_{m_0}}{\partial x}(\nu_m,\tilde x_{m,\nu_m})\frac{d\tilde x_{m,\nu}}{d\nu}(\nu_m)\right|\\
&\geq \left|\frac{\partial y_{m_0}}{\partial \nu}(\nu_m,\tilde x_{m,\nu_m})\right|-
\left|\frac{\partial y_{m_0}}{\partial x}(\nu_m,\tilde x_{m,\nu_m})\frac{d\tilde x_{m,\nu}}{d\nu}(\nu_m)\right|
>C_0
\end{split}
\end{equation}
for some positive constant $C_0$ and all $m$ sufficiently greater than $m_0$.
Since $y_m(\nu,\hat x_m)=\beta_\nu^nh_m(\nu)$ for $\beta_\nu:=\beta_{\mu_0,\nu}$,
\begin{align*}
\frac{\partial y_m}{\partial \nu}(\nu_m,\hat x_m)&=
\beta_{\nu_m}^n\frac{d h_m}{d \nu}(\nu_m)+
n\beta_{\nu_m}^{n-1}\frac{d\beta_\nu}{d\nu}(\nu_m)h_m(\nu_m)\\
&=\beta_{\nu_m}^n\frac{d h_m}{d \nu}(\nu_m)+
\frac{n}{\beta_{\nu_m}}\frac{d\beta_\nu}{d\nu}(\nu_m)y_m(\nu_m,\hat x_{m}).
\end{align*}
Since $\lim_{m\rightarrow \infty}\beta_{\nu_m}=\beta_0>1$ and 
$|y_m(\nu_m,\hat x_{m})|\leq \delta$, the inequality (\ref{e_c}) implies  
$\lim_{m\rightarrow \infty}|{\partial y_m}(\nu_m,\hat x_m)/{\partial \nu}|=\infty$.
This shows (\ref{eqn_1}).
\end{proof}

\begin{proof}[Proof of Theorem \ref{thm_A}]
Propositions \ref{prop1} and \ref{prop2} imply Theorem \ref{thm_A}.  
\end{proof}

\section*{Acknowledgments}
We would like to thank 
Bau-Sen Du and Yi-Chiuan Chen for their support and hospitality, and 
Ming-Chia Li and Mikhail Malkin for their discussions
during the first draft of this paper was written in Academia Sinica of Taiwan.
We also would like to thank the referees for their valuable comments.

\end{document}